\documentclass[pdflatex,sn-mathphys]{sn-jnl}
\newcommand{\reffig}[1]{Fig.\,\ref{#1}}

\newcommand{\reftab}[1]{Table\,\ref{#1}}

\newcommand{\refalg}[1]{Algorithm\,\ref{#1}}

\newcommand{\rmnum}[1]{\romannumeral #1}
\newcommand{\Rmnum}[1]{\expandafter\@slowromancap\romannumeral #1@}

\DeclareMathOperator{\diver}{\nabla \cdot} %

\usepackage{mathrsfs}
\usepackage{amsmath,amssymb,amsfonts} 
\usepackage{cases} 
\usepackage{graphicx,subfigure}
\usepackage{epstopdf} 
\usepackage{bm}
\usepackage{setspace}
\usepackage{multirow}
\usepackage{booktabs} %\toprule, \midrule and \bottomrule
\usepackage{extarrows}
\usepackage{relsize}
\usepackage[section]{placeins} %
\usepackage{xcolor,float} 
\graphicspath{{figs/}}
\usepackage{makecell} % 

\usepackage{algorithm}
\usepackage{algorithmicx}
\numberwithin{algorithm}{section}

\makeatletter
\newcommand*\bigcdot{\mathpalette\bigcdot@{.5}}
\newcommand*\bigcdot@[2]{\mathbin{\vcenter{\hbox{\scalebox{#2}{$\m@th#1\bullet$}}}}}
\makeatother
\jyear{2022}%

\theoremstyle{thmstyleone}%
\theoremstyle{thmstyletwo}%
\newtheorem{remark}{Remark}%

\theoremstyle{thmstylethree}%

\begin{document}
%%=============================================================%%
\title[~]{An Improved Multi-Stage Preconditioner on GPUs for Compositional Reservoir Simulation}


\author*[1]{\fnm{Li} \sur{Zhao}}\email{zhaoli@smail.xtu.edu.cn}

\author[2]{\fnm{Chen-Song} \sur{Zhang}}\email{zhangcs@lsec.cc.ac.cn}

\author[1]{\fnm{Chun-Sheng} \sur{Feng}}\email{spring@xtu.edu.cn}

\author[1]{\fnm{Shi} \sur{Shu}}\email{shushi@xtu.edu.cn}

\affil[1]{\orgdiv{School of Mathematics and Computational Science}, \orgname{Xiangtan University}, \orgaddress{\city{Xiangtan}, \postcode{411105}, \country{P. R. China}}}

\affil[2]{\orgdiv{LSEC \& NCMIS, Academy of Mathematics and Systems Science, CAS}, \orgaddress{\city{Beijing}, \postcode{100190}, \country{P. R. China}}}

%\affil[3]{\orgdiv{Department}, \orgname{Organization}, \orgaddress{\street{Street}, \city{City}, \postcode{610101}, \state{State}, \country{Country}}}


\abstract{
The compositional model is often used to describe multicomponent multiphase porous media flows in the petroleum industry. The fully implicit method with strong stability and weak constraints on time-step sizes is commonly used in the mainstream commercial reservoir simulators. In this paper, we develop an efficient multi-stage preconditioner for the fully implicit compositional flow simulation. The method employs an adaptive setup phase to improve the parallel efficiency on GPUs. Furthermore, a multi-color Gauss--Seidel algorithm based on the adjacency matrix is applied in the algebraic multigrid methods for the pressure part. 
Numerical results demonstrate that the proposed algorithm achieves good parallel speedup while yields the same convergence behavior as the corresponding sequential version. 
}




\maketitle

\section{Introduction}\label{sec:1}
The compositional model, which allows the fluids to be composed of various material components, is a widely-used mathematical model for describing multiphase flows in porous media~\cite{Aziz1979PetroleumRS,ChenHuanMa}. The compositional model is an extension of the blackoil model~\cite{Peaceman1977}, which is formed by multiple coupled nonlinear partial differential equations. Some complicated oil displacement technologies can be accurately simulated based on the compositional model, such as polymer flooding, surfactant and alkali oil displacement agents, and miscible flooding.

Numerical methods for compositional numerical simulation are abundant; to name a few, IMplicit Pressure Explicit Concentrations (IMPEC) method~\cite{Fussell1979}, Fully Implicit Method (FIM)~\cite{Coats1980}, IMplicit Pressure / SATuration and explicit concentrations (IMPSAT)~\cite{quandalle1989implicit}, and Adaptive Implicit Method (AIM)~\cite {AIM}.
In the IMPEC method, pressure is implicit and other variables are explicit. One of its advantages is that no need to solve coupled linear algebra systems, but its time stepsize is constrained by the Courant--Friedrichs--Lewy (CFL)~\cite{Courant1928} condition. The FIM method, on the contrary, is unconditionally stable with respect to time stepsizes because all variables are handled implicitly. The IMPSAT method is a combination of the IMPEC and FIM methods, where pressure and saturation variables are handled implicitly and the molar fractions are calculated explicitly. The AIM method is also a compromise between the IMPEC and FIM methods and it also yields Jacobian algebraic systems that are easier to solve than FIM.

The FIM method possesses the characteristics of good stability and is widely used in the commercial simulators. However, a coupled nonlinear system of equations need to be solved at each time step of FIM. Such systems are usually linearized using the Newton-type methods, which requires solving a coupled Jacobian linear algebra system during each iteration. In the numerical simulation, the solution of these systems is the main computational cost~\cite{Chensong2022}. Therefore, efficient linear solvers are crucial for improving the efficiency of fully implicit reservoir simulation, especially for large-scale three-dimensional problems.

The linear solution methods generally consist of a setup phase (SETUP) and a solve phase (SOLVE). Iterative methods are widely used in petroleum reservoir simulation due to their low memory overhead and good parallel scalability. More specifically, the GMRES and BiCGstab methods~\cite{2003Iterative} are exploited to solve the nonsymmetric systems arisen from the fully implicit discretization of reservoir models. The preconditioning techniques are crucial to speed up the convergence of iterative methods~\cite{Chensong2022}. 
For large-scale reservoir simulation, multi-stage preconditioners are very competitive. The classical Constrained Pressure Residual (CPR)~\cite{10.2118/96809-MS,2017Numerical,IGE1983,10.2118/13536-MS} approach is a well-known two-stage preconditioner. It utilizes the Algebraic MultiGrid (AMG)~\cite{1984Algebraic,falgout2006an} method to approximate the inverse of the pressure matrix in the first stage and  the Incomplete LU (ILU) factorization~\cite{10.2118/12262-MS} to smooth the overall reservoir matrix in the second stage. The Multi-Stage Preconditioner (MSP)~\cite{10.2118/118722-MS, XiaoZheHU2013,10.2118/105832-MS} is a generalization of the CPR method, which is also widely used in petroleum reservoir problems.

Parallel computing is an important approach to improve the speed of simulation and a lot of attention has been paid to developing efficient parallel algorithms. A Graphics Processing Unit (GPU) with thousands of cores is a parallel accelerator that are designed to handle images and graphics originally. Due to its high float-point performance and memory bandwidth~\cite{NVIDIA2022}, it has great potential in petroleum reservoir simulation. In recent years, research on GPU parallel algorithms has been developed in~\cite{2010High,chen2014gpu,YangBo2016,kang2018parallel,Manea2019massivelyGPU,Middya2021GPU,Esler2022GPU} and the references therein. For example, Chen et al.~\cite{chen2014gpu, YangBo2016} developed a hybrid sparse matrix storage format and the corresponding sparse matrix-vector multiplication (SpMV) kernel. Kang et al.~\cite{kang2018parallel} developed a parallel nonlinear solver based on OpenACC~\cite{OpenACC2022} using the domain decomposition method to achieve load balancing. Finally, Manea et al.~\cite{Manea2019massivelyGPU} studied a parallel algebraic multiscale solver on GPU architectures to improve the solution efficiency of the pressure equation.

In this paper, we focus on a GPU-based parallel linear solver for compositional models, which is an extension of the recent work for the blackoil model~\cite{LiZhao2022}. To the best of our knowledge, few numerical tests have been done for the algebraic systems arising from the fully implicit discretization of the compositional model in the literature. For such systems, we develop a parallel multi-stage preconditioner. To begin with, we propose a multi-stage preconditioner with an adaptive SETUP, denoted as ASMSP. The proposed method can significantly reduce the number of SETUP calls, so as to reduce the computational overhead of SETUP and improve parallel efficiency. 
Furthermore, we investigate a multi-color Gauss--Seidel (GS) algorithm based on algebraic grouping for the smoothing operator in the AMG methods. The proposed algorithm can produce the same convergence behavior as the corresponding single-threaded algorithm and is suitable for unstructured grids.
Finally, the proposed methods are integrated into our open-source simulator OpenCAEPoro~\cite{OpenCAEPoro} for multicomponent multiphase flows in porous media.

The rest of the paper is organized as follows. Section \ref{sec:2} briefly introduces the compositional model and its fully implicit discretization. In Section \ref{sec:3}, an MSP preconditioner with adaptive SETUP is developed. In Section \ref{sec:4}, the parallel implementation of multi-color GS based on the adjacency matrix is proposed. In Section \ref{sec:5}, numerical experiments are performed to evaluate the convergence and parallel speedup of the proposed method. Section \ref{sec:6} summarizes the work of this paper.

\section{Mathematical model and discretization}\label{sec:2}
\subsection{The compositional model}\label{sec:2-1}
In this paper, we consider the isothermal multicomponent compositional model~\cite{Aziz1979PetroleumRS,ChenHuanMa} containing $ n_c $ components (hydrocarbon and water) and $ n_p $ phases (including at least the water phase). The mass conservation equation for the component $ i $  reads
\begin{equation}\label{eq:MassConservation}
	\frac{\partial}{\partial t}\Big(\phi \sum\limits_{j=1}^{n_p}x_{ij} \xi_j S_j \Big) + \diver \Big( \sum\limits_{j=1}^{n_p}x_{ij} \xi_j \bm{u}_j \Big) = Q_{i},~i=1,2,\ldots,n_c,	
\end{equation}
where $ \phi $ is the porosity of the rock, $ x_{ij} $ is the molar fraction (dimensionless) of component $i$ in phase $j$, $ \xi_j $ is the molar concentration of phase $ j $, $ S_{j}$ is the saturation of phase $ j $, $ \bm{u}_{j} $ is the velocity of phase $ j $, and $ Q_{i} $ is source/sink terms of component $i$.

Assume that the pore volume of the porous media is filled with the fluid, the volume balance equation is then
\begin{equation}\label{eq:VolumeBalance}
	V = V_{\rm{pore}} := \phi V_{\rm{bulk}},
\end{equation}
where $ V $ is the fluid volume, $ V_{\rm{pore}} $ is the pore volume, and $ V_{\rm{bulk}} $ is the bulk volume.

Assume that the phase $ j $ fluid in porous media satisfies the Darcy's law:
\begin{equation}\label{eq:Darcy}
	\bm{u}_{j}=-\frac{\kappa \kappa_{r j}}{\mu_{j}}\left(\nabla P_{j}-\rho_{j} \mathfrak{g} \nabla z\right),~j=1,2,\ldots,n_p,
\end{equation}
where $ \kappa $ is the absolute permeability, $ \kappa_{r j} $ is the relative permeability of phase $j $, $ \mu_{j} $ is the viscosity coefficient of phase $ j $, $ P_{j }$ is the pressure of phase $ j $, $ \rho_{j} $ is the density of phase $ j $, $ \mathfrak{g} $ is the gravity acceleration, and $ z $ is the depth.

Moreover, the variables $S_{j}$, $ x_{ij} $ and $P_{j}$ in the Eqs.~\eqref{eq:MassConservation}-\eqref{eq:Darcy} satisfy the following constitutive relations:
\begin{itemize}
	\item Saturation constraint equation:
	\begin{equation}\label{eq:saturation}
		\sum_{j=1}^{n_p} S_{j} = 1,
	\end{equation}
	
	\item Molar fraction constraint equation:
	\begin{equation}\label{eq:molar-fraction}
		\sum_{i=1}^{n_c} x_{ij} = 1,~j=1,2,\ldots,n_p,
	\end{equation}
	
	\item Capillary pressure equation:
	\begin{equation}\label{eq:capillary}
		P_j = P - P_{cj}(S_j),~j=1,2,\ldots,n_p,
		%P = P_j + P_{cj}(S_j) \approx P_j,~j=1,2,\ldots,n_p,
	\end{equation}
	where $ P $ is the reference pressure, and $ P_{cj}(S_j) $ is the capillary pressure between the reference phase and phase $ j $, which will be ignored in the rest of this paper.
\end{itemize}

\subsection{Discretization method}
In this section, we first simplify the compositional model, then describe the choice of main equations and primary variables, and finally present the FIM discretization method. 

\subsubsection{The choices of primary variables}
To begin with, we introduce the overall molar concentration $ N_i $ and molar flux $ \bm{F}_i $ of the component $ i $, which are defined as
\begin{equation}\label{eq:Ni-Fi}
	N_i = \phi \sum\limits_{j=1}^{n_p}x_{ij} \xi_j S_j,
\end{equation}
\begin{equation}\label{Fi}
	\bm{F}_i = \sum\limits_{j=1}^{n_p}x_{ij} \xi_j \bm{u}_j.
\end{equation}
Eq.~\eqref{eq:MassConservation} can be simplified to
\begin{equation}\label{eq:SimMassConservation}
	\frac{\partial}{\partial t} N_i + \diver \bm{F}_i = Q_{i},~i=1,2,\ldots,n_c,	
\end{equation}

In this paper, we choose the mass conservation equations~\eqref{eq:SimMassConservation} and the volume balance equation~\eqref{eq:VolumeBalance} as the main equations, and there are $n_c + 1$ equations in total. The reference pressure $ P $ and the overall molar concentration $ N_i~(i=1,2,\ldots,n_c) $ are used as the primary variables of the discrete method, and there are $n_c + 1$ variables in total. After solving the primary variables, the fluid volume state function $V(P,N_1,\ldots,N_{n_c})$ can be obtained by the equation of state and flash calculation, see~\cite{Peng1977,ChenHuanMa} for more details.

\subsubsection{Finite volume method and backward Euler method}
The finite volume method (FVM) \cite{LeVeque2002} features simplicity, conservation, adaptivity to complex geometric regions, and monotonicity; it is a commonly used discretization method in the petroleum industry. In this paper, the spatial domain is discretized by using FVM. Suppose that the spatial domain $ \Omega \subset \mathbb{R}^3 $ (an open set) is discretized into $ m $ elements set $ \{ \tau_k \}_{k=1}^{m} $ (the shape of the element is not considered here), which satisfies $ \cup_{k=1}^{m} {\overline \tau_k} = \overline{\Omega}$ and $ \tau_k \cap \tau_\ell = \varnothing,~k \neq \ell,~k,\ell =1,\ldots,m$.

\begin{equation}\label{eq:IntMC}
	\int_{\tau_k} \frac{\partial}{\partial t} N_i dV + \int_{\tau_k} \diver \bm{F}_i dV = \int_{\tau_k} Q_{i} dV,~i=1,2,\ldots,n_c,
\end{equation}
using the divergence theorem, we can get
\begin{equation}\label{eq:IntMC-div}
	\int_{\tau_k} \frac{\partial}{\partial t} N_i dV + \int_{S_k} \bm{F}_i \cdot \bm{n} dS = \int_{\tau_k} Q_{i} dV,~i=1,2,\ldots,n_c,
\end{equation}
where $ S_k := \partial \tau_k$ is all the surfaces set of element $ \tau_k $ and $ \bm{n} $ is the outer unit normal vector to $ S_k $.

Therefore, the discrete equation on element $ \tau_k $ can be written as
\begin{equation}\label{eq:semi-discrete-eqs}
	\frac{\partial}{\partial t} N_{i,k} + \sum_{s \in S_k} \bm{F}_{i,s} = Q_{i,k},~i=1,2,\ldots,n_c,
\end{equation}
where the discrete flux $ \bm{F}_{i,s} $ can be defined in various ways (e.g., \cite{Aavatsmark2002MPFA,Aavatsmark2008MPFA}), we consider the following form
\begin{equation}\label{eq:flux}
	\bm{F}_{i,s} = \left\{ \frac{L \kappa}{d} \right\}_s \sum\limits_{j=1}^{n_p} 
	\left( \left\{ \frac{\kappa_{rj}}{\mu_{j}} \right\}_s \big\{ x_{ij} \xi_j \big\}_s 
	\delta_s (P+P_{cj}-\rho_{j} \mathfrak{g} z ) \right), 
\end{equation}
Here, $L$ and $d$ denote the size of the interface and the distance between two adjacent elements, respectively. $ \delta_s $ denotes the difference between the values on the two adjacent elements. Because the primary variables of the discrete equations are defined at the center of element, the value $\{\cdot\}_s$ of the physical quantities on the interface $s$ needs to be approximated by some suitable methods. Typically, the harmonic mean value or the upstream weighted value of the physical quantities in the two elements of the shared interface $s$ is taken; see \cite{ChenHuanMa} for more details.

For the semi-discrete equation \eqref{eq:semi-discrete-eqs}, the time derivative term is discretized by using the backward Euler method, and the superscripts $n$ and $n + 1$ denote the time $t_n$ and $t_{n+1}$, respectively. The fully discrete volume balance equation and mass conservation equations are, respectively.
\begin{equation}\label{eq:DiscreteVolumeBalance}
	V^{n+1} - V_{\rm{pore}}^{n+1} = 0,
\end{equation}
\begin{equation}\label{eq:fully-discrete-eqs}
	\frac{N_{i,k}^{n+1} - N_{i,k}^{n} }{\Delta t} + \sum_{s \in S_k} \bm{F}_{i,s}^{\star} = Q_{i,k}^{\star},~i=1,2,\ldots,n_c,
\end{equation}
where the source/sink term is simplified into a known function; but in practical problems, it is related to the production mode of the well in the oil field, and is also strongly coupled with the primary variables. Since the focus of this paper is not how to handle the well equations, we do not describe it in detail.

Note that the fully discrete equations \eqref{eq:DiscreteVolumeBalance} and \eqref{eq:fully-discrete-eqs} are nonlinear, and the terms $ \bm{F}_{i,s}^{\star} $ and $ Q_{i,k}^{\star} $ are subject to be specified. Below, we will give their expressions.

\subsubsection{Fully implicit method}
The FIM scheme is currently commonly used in mainstream commercial reservoir simulators. This is because the scheme has the characteristics of strong stability and weak constraint on the timestep sizes. These characteristics highlight the advantages of the FIM, especially when the nonlinearity of the models is relatively strong. 

When both $ \bm{F}_{i,s}^{\star} $ and $ Q_{i,k}^{\star} $ in Eq.~\eqref{eq:fully-discrete-eqs} take the value of $t_{n+1}$ time, the fully implicit discrete equations are
\begin{equation}\label{eq:FIM-nonlinear}
	\frac{N_{i,k}^{n+1} - N_{i,k}^{n} }{\Delta t} + \sum_{s \in S_k} \bm{F}_{i,s}^{n+1} = Q_{i,k}^{n+1},~i=1,2,\ldots,n_c,
\end{equation}

Owing to the implicit solution for $n_c + 1$ primary variables, Eqs.~\eqref{eq:DiscreteVolumeBalance} and \eqref{eq:FIM-nonlinear} are strongly coupled nonlinear systems of equations that need to be linearized. In this work, we exploit the well-known Newton's method to linearize Eqs.~\eqref{eq:DiscreteVolumeBalance} and \eqref{eq:FIM-nonlinear}. The Jacobian equation for increments $\delta P,\delta N_1,\ldots,\delta N_{n_c}  $ can be written as
\begin{equation}\label{eq:FIM-Newton-linear1}
	\frac{1}{\Delta t} \alpha_P \delta P -\frac{1}{\Delta t} \sum_{i=1}^{n_c} \alpha_i \delta N_i = r_P,
\end{equation}
\begin{equation}\label{eq:FIM-Newton-linear2}
	\frac{1}{\Delta t} \delta N_i - \diver (\beta_i \nabla \delta P)
	- \diver (\beta_{iP} \nabla \delta P) - \sum_{k=1}^{n_c} \diver (\beta_{ik} \nabla \delta N_k)= r_i,~i=1,\ldots,n_c,
\end{equation}
where the coefficients $ \alpha_P $, $ \alpha_i $, $ \beta_i $, $ \beta_{iP} $, and $ \beta_{ik} $ are obtained by partial derivation of the model coefficients with respect to $P$ or $N_i$; see \cite{Qiao2015} for details.

\subsubsection{Discrete system}
After discretization, the coupled nonlinear algebraic equations are obtained. 
Such equations are linearized by adopting the Newton method to form the sparse Jacobian system $ Ax = b $ of the reservoir equation with implicit wells, namely:
\begin{equation}\label{eq:Jacobian-system}
	\left( {\begin{array}{*{20}{c}}
			{{{{A}}_{RR}}}&{{A_{RW}}}\\
			{{A_{WR}}}&{{A_{WW}}}
	\end{array}} \right)
	\left( {\begin{array}{*{20}{c}}
			{{x_R}}\\
			{{x_W}}
	\end{array}} \right) 
	= \left( {\begin{array}{*{20}{c}}
			{{b_R}}\\
			{{b_W}}
	\end{array}} \right),
\end{equation}
where $ A_{RR} $ and $ A_{RW} $ are the derivatives of the reservoir equations for reservoir variables and well variables, respectively; $ A_{WR} $ and $ A_{WW} $ are the derivatives of the well equations for reservoir variables and well variables, respectively; $ x_R $ and $ x_W $ are reservoir and bottom-hole flowing pressure variables, respectively; and $ b_R $ and $ b_W $ are the right-hand side vectors that correspond to the reservoir fields and the implicit wells, respectively.

The subsystem corresponding to the reservoir equations in the discrete system \eqref{eq:Jacobian-system} is $ A_{RR}x_R=b_R $; that is,
\begin{equation}\label{eq:Jacobian-system-reservoir}
	\left( {\begin{array}{*{20}{c}}
			A_{PP}& A_{PN_1} & A_{PN_2} & \cdots &A_{PN_{n_c}} \\
			A_{N_1P}& A_{N_1N_1} & A_{N_1N_2} & \cdots &A_{N_1 N_{n_c}} \\
			A_{N_2P}& A_{N_2N_1} & A_{N_2N_2} & \cdots &A_{N_2 N_{n_c}} \\
			\vdots & \vdots & \vdots & \vdots & \vdots \\
			A_{N_{n_c}P}& A_{N_{n_c}N_1} & A_{N_{n_c} N_2} & \cdots &A_{N_{n_c} N_{n_c}} \\
	\end{array}} \right)%,
	\left( {\begin{array}{*{20}{c}}
			x_{P}   \\
			x_{1} \\
			x_{2} \\
			\vdots  \\
			x_{n_c} \\
	\end{array}} \right)%, 
	= 
	\left( {\begin{array}{*{20}{c}}
			b_{P}   \\
			b_{1} \\
			b_{2} \\
			\vdots  \\
			b_{n_c} \\
	\end{array}} \right),
\end{equation}
where $ P $ is reference pressure and $ N_i~(i=1,2,\ldots,n_c) $ are the overall molar concentration.

\begin{remark}
	\rm{For convenience, we do not describe how to deal with the well equations.}
\end{remark}

\section{Parallel multi-stage preconditioners with adaptive SETUP}\label{sec:3}

In our compositional model, the primary variables consist of reference pressure $ P $ and overall molar concentration $ N_i~(i=1,2,\ldots,n_c) $. These variables possess different mathematical properties, such as the parabolicity of the pressure equation and hyperbolicity of the concentration equations. These properties provide a theoretical basis for the construction of multiplicative subspace correction methods~\cite{Xu1992MSC, Xu1996ASP}; see~\cite{Chensong2022} for a recent review. 

\subsection{Multi-stage preconditioner}

We first define two transfer operators of the reservoir matrix, suppose $ \varPi_N: \mathcal{V}_N \rightarrow \mathcal{V} $ and $ \varPi_P: \mathcal{V}_P \rightarrow \mathcal{V} $, where $ \mathcal{V}_N $ and $ \mathcal{V}_P $ are the overall molar concentration and pressure variables space, respectively, and $\mathcal{V}$ is the variables space of the whole reservoir. Then, the multiplicative multi-stage preconditioner $ B $~\cite{10.2118/118722-MS, XiaoZheHU2013,10.2118/105832-MS} is defined as
\begin{equation}\label{eq:MSP}
	I - B A = (I-RA)(I- \varPi_P B_P \varPi_P^T A)(I- \varPi_N B_N \varPi_N^T A).
\end{equation}
where the relaxation operator $ R $ employs the Block ILU (BILU) method, $B_P$ and $B_N$ are solved by the AMG and Block GS (BGS) methods, respectively.

Suppose that the mathematical behavior of preconditioner $ B $ acting on a known vector $ g $ is 
\begin{equation}\label{eq:MSP-behavior}
	w=B g,
\end{equation}
The corresponding multi-stage preconditioning algorithm \cite{2014A2} is shown in the \refalg{alg:MSP}.

\begin{algorithm}[H]\label{alg:MSP}
	\caption{MSP preconditioning method}
	\begin{algorithmic}[1]
		\Require $A, g, w, \varPi_N, \varPi_P$;
		\Ensure $w = Bg$.
		\State $ r = g-Aw$;
		\State $ w = w + \varPi_N B_N \varPi_N^{T}r $; 
		\State $ r = g-Aw$;
		\State $ w = w + \varPi_P B_P \varPi_P^{T}r $; 
		\State $ r = g-Aw$;
		\State $ w = w + Rr $.
	\end{algorithmic}
\end{algorithm}

\subsection{MSP with adaptive SETUP}
\begin{figure}[H]
	\centering
	\includegraphics[width=0.6\linewidth]{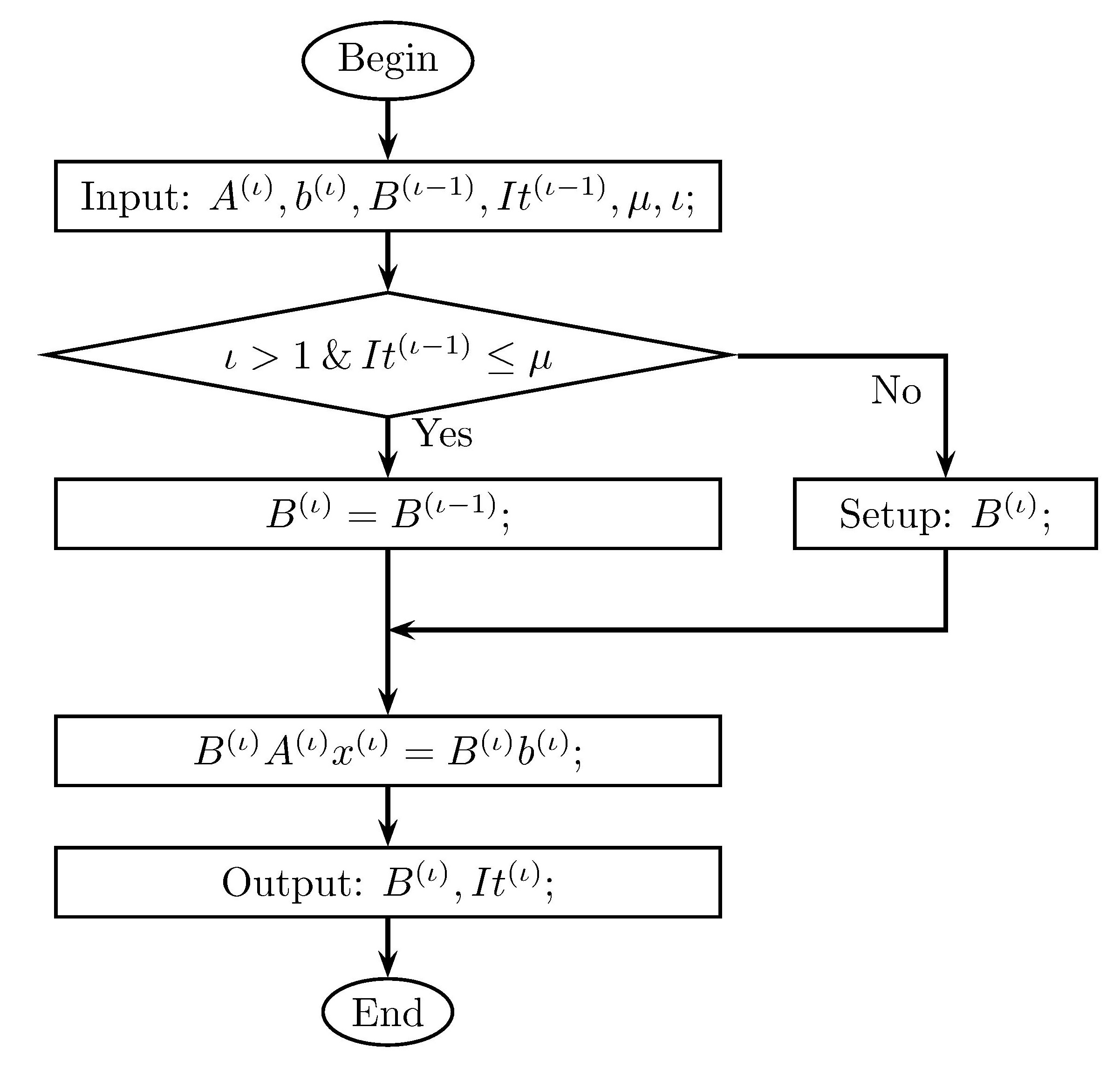}
	\caption{The algorithm flow chart of adaptive SETUP.}
	\label{fig:sharesetup}
\end{figure}

In this subsection, we utilize an adaptive SETUP strategy for the MSP preconditioner to improve its parallel efficiency, denoted as ASMSP. As mentioned earlier, this strategy has been employged for the blackoil model by Zhao et al.~\cite{LiZhao2022}. Assume that the solution objective is $A^{(\iota)}x^{(\iota)}=b^{(\iota)}$. It is worth mentioning that superscript $\iota$ is the number of Newton iterations. This is because the reservoir model is a nonlinear system of equations, which is linearized using Newton's method here (see~\cite{ChenHuanMa} for more details). \reffig{fig:sharesetup} presents the algorithm flow chart of adaptive SETUP. 

In \reffig{fig:sharesetup}, the main difference from the standard methods is that the preconditioner $B^{(\iota)}$ is yielded by an adaptive strategy. This strategy can be divided into the following two cases:
\begin{enumerate}
	\item [(1)] The preconditioner $ B^{(\iota)} $ inherits the information from the previous preconditioner $ B^{(\iota-1)} $. A natural approach is to use the number of iterations $ It^{(\iota-1)} $, required by solving the previous Jacobian system $ A^{(\iota-1)} x^{(\iota-1)} = b^{(\iota-1 )}$. We introduce a threshold $\mu$ (a non-negative integer); if $ It^{(\iota-1)} \leq \mu$, the previous preconditioner $ B^{(\iota-1)} $ is used as the preconditioner $ B^{(\iota)} $.

	\item [(2)] The preconditioner $ B^{(\iota)} $ is regenerated; if $ \iota=1 $ or $ It^{(\iota-1)} > \mu$, the preconditioner $ B^{(\iota)} $ is generated by calling \refalg{alg:MSP}.
\end{enumerate}

\begin{remark}\label{rmk:Regenerate}
	\rm{If the sizes of $ A^{(\iota-1)} $ and $ A^{(\iota)} $ are not the same, the preconditioner $ B^{(\iota)} $ must be regenerated for sure.}
\end{remark}

Below, we illustrate the rationale for this approach. The number of iterations can evaluate the quality of a preconditioner. More iterations indicate a poor preconditioner, and fewer iterations indicate a good preconditioner. Especially when the number of iterations is 1, the preconditioner is the exact inverse of the matrix. The $ It^{(\iota-1)} \leq \mu $ indicates $ B^{(\iota-1)} $ is an effective preconditioner for $ A^{(\iota-1)} x^{(\iota-1)} = b^ {(\iota-1)}$. In addition, the structure of these matrices is very similar during Newton's iteration, and the preconditioner does not need to approximate the inverse of the matrix exactly. So the preconditioner $B^{(\iota-1)}$ can also be applied to the Jacobian system $A^{(\iota)} x^{(\iota)} = b^{(\iota)}$. The proposed method can improve the parallel performance of the solver by reducing the number of SETUP calls and reducing the proportion of low parallel speedup in the solver.

Finally, we discuss the impact of the threshold $\mu $ on performance. If $\mu $ is too small, the number of SETUP calls will not be significantly reduced, which will not significantly improve the performance of the solver. In particular, ASMSP degenerates to standard MSP when $\mu =0$. Conversely, if $\mu $ is too large, too few SETUP calls can also affect the performance, due to the dramatic increase in the number of iterations. Usually, a suitable $\mu $ is determined by numerical experiments.

\section{A multi-color GS based on adjacency matrix}\label{sec:4}
It is well-known that the GS algorithm, compared to the Jacobi algorithm, exploits most updated values in the iterative process. Therefore, the GS algorithm brings a better convergence. However, it is essentially sequential and cannot be easily parallelized. A popular red-black GS (also known as multi-color GS) parallel algorithm has attracted a lot of attention~\cite{2003Iterative}. Unfortunately, the algorithm is designed based on structured grids and is not compatible with unstructured grids. 

A hybrid approach that combines the Jacobi and GS methods can be applied, but its convergence rate also deteriorates with respect to higher parallelism. In order to overcome the limitations of the traditional red-black GS algorithm, a multi-color GS algorithm based on the coefficient matrix of strong connections has been proposed and analyzed in~\cite{LiZhao2022}. This paper proposes a multi-color GS algorithm from the algebraic point of view. The proposed method yields the same convergence behavior as the corresponding single-threaded algorithm; moreover, it obtains good parallel performance when using a lot of threads on GPUs.

\subsection{Adjacency graph and algorithm principles}\label{sec:alg-target}
The notion of an adjacency graph needs to be introduced to implement the multi-color GS algorithm algebraically. An adjacency graph corresponds to a sparse matrix, which reflects the nonzero pattern of the matrix, i.e., the nonzero entries of the matrix reflect the connectivity relationship between the vertices in the adjacency graph.

We develop a multi-color GS algorithm that can be applied to symmetric and nonsymmetric matrices. For simplicity, assuming that the sparse matrix $A \in \mathbb{R}^{n\times n} $ is a symmetric matrix. Let $G_{A}(V, E)$ be the adjacency graph corresponding to the matrix $A=\left(a_{ij}\right)_{n \times n}$. Here $V=\left\{ v_{1}, v_{2}, \ldots, v_{n}\right\}$ and $E =\left\{\left(v_{i}, v_{j}\right): \forall \, i \neq j, a_{ij} \neq 0 \right\}$ are the vertices and edges sets, respectively. It is easy to know that each nonzero entry $a_{ij}$ on the off-diagonal of $A$ corresponds to an edge $(v_{i}, v_{j})$.

Here, we give the principles for designing a multi-color GS algorithm in this paper:
\begin{enumerate}
	\item [(\rmnum{1})] The vertices set $V$ is divided into $g$ subsets: $V = V_{1} \cup V_{2} \cdots \cup V_{g}$;
	\item [(\rmnum{2})] Any two subsets are disjoint: $V_{i} \cap V_{j}=\varnothing,~i \neq j,~i,j = 1,2, \ldots, g$;
	\item [(\rmnum{3})] Vertices in any subset are not connected by edges: $a_{ij} =a_{ji}=0, ~\forall \, v_{i}, v_{j} \in V_{\ell},~\ell=1,2,\ldots,g$; 
	\item [(\rmnum{4})] The number of subsets, $g$, should be as small as possible.
\end{enumerate}
It is obvious that the difficulty of grouping and parallel granularity increase as the number of groups $g$ decreases. In fact, the red-black GS algorithm satisfies the conditions for $g=2$. In particular, when $g=n$, it degenerates to the classic GS algorithm.

\subsection{Multi-color GS algorithm}
We define the adjacency matrix $S$ corresponding to the adjacency graph $G_{A}(V, E)$. Its diagonal entries are zero and off-diagonal entries are
\begin{equation}\label{eq:adjacency-matrix}
	S_{i j}=\left\{\begin{array}{ll}
		1, & \text{if}~a_{i j} \neq 0, \\
		0, & \text{if}~a_{i j} = 0,
	\end{array} \quad \forall \, i, j=1,2, \ldots, n,~i \neq j, \right.
\end{equation}
where $S_{ij}=1$ denotes an adjacent edge between $v_{i}$ and $v_{j}$, and $S_{ij}=0$ denotes no adjacent edge between $v_{i}$ and $v_{j}$.
To describe the multi-color GS algorithm, some notations are introduced in \reftab{tab:notations}.

\begin{table}[h]
\centering	
\renewcommand\arraystretch{1.5} 
\caption{The definitions and explanations of some notations.}\label{tab:notations}
\begin{tabular}{c|c}
\hline
Notation & Definition and explanation \\
\hline
\multirow{2}{*}{$S_{i}$}
&\makecell[l]{$ S_i = \left\{ j: S_{ij} \neq 0,j=1,2, \ldots, n \right\} $,} \\
&\makecell[l]{$S_{i}$ denotes the vertex set that is connected to the vertex $v_{i}$.} \\
\hline
\multirow{2}{*}{$\overline{S}_{i}$}
&\makecell[l]{$ \overline{S}_i = \big\{ j: j\in S_{i} \cup \{i\}~\text{and}~ \text{color of}~j~\text{is undetermined} \big\} $,} \\
&\makecell[l]{$\overline{S}_{i}$ denotes the vertex set that is connected to the vertex $v_{i}$ (including $v_{i}$) and \\ whose colors are undetermined.} \\ 
\hline
\rule{0pt}{15pt} 
\multirow{3}{*}{$\widehat{S}_{i}$}
&\makecell[l]{$ \widehat{S}_i = \big\{ j: j\in W_{i}~\text{and}~\text{color of}~j~\text{is undetermined} \big\}$, where \\$W_i := \big\{ j: \forall\,k \in S_i, j \in S_k / (S_i \cup \{i\}) \big\} $,} \\
&\makecell[l]{$\widehat{S}_{i}$ denotes the vertex set that is the next connected to the vertex $v_{i}$ \\ (the vertices on ``the second circle") and whose colors are undetermined.} \\ 
\hline
\multirow{2}{*}{$\lvert S_{i} \rvert $}
&\makecell[l]{$\lvert S_{i} \rvert = \sum_{j \in S_i} 1 $,} \\
&\makecell[l]{The cardinality $\lvert S_{i} \rvert $ denotes the number of entries in the set $ S_{i} $.} \\
\hline
\end{tabular}
\end{table}

Below, we first present a (greedy) splitting algorithm for the set of vertices $ V $ based on the adjacency matrix $ S $, denoted as VerticesSplitting (see \refalg{alg:strong-vsplit}). Then we give the vertices grouping algorithm of matrix $ A $, denoted as VerticesGrouping (see \refalg{alg:v-split}). Finally, a multi-color parallel GS method is presented in~\refalg{alg:GS-parallel}, denoted as PGS-MC.

\begin{algorithm}
%{\small{	
	\caption{VerticesSplitting method} 
	\label{alg:strong-vsplit}
	\begin{algorithmic}[1]
%		\State \textbf{Input:} $V,~S$;
		\Require $V,~S$;
		\Ensure $ W,~\overline{W}$.
		
		\State Let $ W = \varnothing,~\overline{W} = \varnothing,~ \widehat{W} = \varnothing $;
		
		\While{$ V \neq \varnothing $}
		\If{$ \widehat{W} \neq \varnothing $}
		\State Any take $ v_i \in \widehat{W}$ and $ \lvert S_i \rvert \geq  \lvert S_j \rvert$, $ \forall \, v_i, v_j \in \widehat{W} $;
		\Else
		\State Any take $ v_i \in V$ and $ \lvert S_i \rvert \geq  \lvert S_j \rvert$, $ \forall \, v_i, v_j \in V $;
		\EndIf
		\If{$ v_i $ is not connected to any vertices in $ W $ (i.e., $S_{ij} = 0,\forall \, j \in W$)}	
		\State $ W = W \cup v_i,~V = V/v_i $;
		
		\If{$ v_i \in \widehat{W} $}
		\State $ \widehat{W} = \widehat{W}/v_i $;
		\EndIf	
		
		\State $ \overline{W} = \overline{W} \cup S_i $,~$ V \leftarrow V / S_i $,~$ \widehat{W} = \widehat{W} \cup \widehat{S}_i $;		
%		\State $ V \leftarrow V / S_i $;		
%		\State $ \widehat{W} = \widehat{W} \cup \widehat{S}_i $;
		\Else
		\State $ \overline{W} = \overline{W} \cup v_i $,~$ V \leftarrow V / v_i $;
%		\State $ V \leftarrow V / v_i $;
		
		\If{$ v_i \in \widehat{W} $}
		\State $ \widehat{W} = \widehat{W}/v_i $;
		\EndIf
		
		\EndIf	
		\EndWhile
%		\State  \textbf{Output: } $ W,~\overline{W}$. 
	\end{algorithmic} 
%}}
\end{algorithm} 

\begin{algorithm}
	\caption{VerticesGrouping method}
	\label{alg:v-split}
	\begin{algorithmic}[1]
%		\State \textbf{Input:} $V,~S$;
		\Require $V,~S$;
		\Ensure $ \{ V_1,V_2,\ldots,V_g \}$.
		\State Set $g = 0$; 
		\While{$ V \neq \varnothing$}
		\State $g = g + 1$;
		\State Get $V_g$ and $\overline{V}_g$ by calling \refalg{alg:strong-vsplit};
		\State Let $ V = \overline{V}_g$;
		\EndWhile
%		\State  \textbf{Output: } $ V_\ell~(\ell=1,\ldots,c)$. 
	\end{algorithmic} 
\end{algorithm} 

In our PGS-MC method, it is worth noting that $G_{A}(V, E)$ can be split into $ g $ subgraphs $G_{A_\ell}(V_\ell, E_\ell)$ by calling~\refalg{alg:v-split}, and the adjacency matrice corresponding to these subgraphs are $ S_\ell$ ($\ell=1,\ldots,g$). It is easy to know that the submatrix $A_\ell$ (corresponding to the subgraph $G_{A_\ell}(V_\ell, E_\ell)$) is a diagonal matrix. In summary, the proposed method starts from the adjacency matrix of the coefficient matrix and designs a vertices grouping algorithm. This method can run in parallel within the same group. Moreover, the proposed method can be applied to the AMG methods, and the numerical experiments in the next section also present its parallel performance.

\begin{algorithm}
	\caption{PGS-MC method}
	\label{alg:GS-parallel}
	\begin{algorithmic}[1]
		%		\State \textbf{Input:} $ A, x, b, \theta $;
		\Require $ A, x, b$;
		\Ensure $ x $. 
		\State Create the vertices set $ V $ and the adjacency matrix $ S $ using the matrix $ A $ and the formula \eqref{eq:adjacency-matrix};
		\State Generate independent vertices subset $V_{\ell}$ ($\ell=1,\ldots,g$) by calling~\refalg{alg:v-split};
		\State Use $V_{\ell}$ to split the matrix $ A $ into submatrix $ A_\ell$;
		\For{$\ell=1, \ldots, g$}
		\State Call classical GS algorithm in parallel for submatrix $A_\ell$;
%		\State Parallel call the classic GS algorithm of submatrix $A_\ell$;
		\EndFor
		%		\State  \textbf{Output: } $ x $. 
	\end{algorithmic} 
\end{algorithm}

\section{Numerical experiments}\label{sec:5}
In this section, benchmark problems based on SPE1, SPE5, and SPE10~\cite{Odeh1981-SPE1,Killough1987SPE5,2001TenthSPE} are considered to demonstrate the performance of the proposed methods. Note that the SPE1 and SPE10 problems can be solved using a blackoil framework as well; but we solve them using the compositional simulator OpenCAEPoro. 

In the ASMSP-GMRES method, the Unsmoothed Aggregation AMG (UA-AMG) method is used to approximate the inverse of the pressure coefficient matrix, where the aggregation strategy is the so-called nonsymmetric pairwise matching aggregation (NPAIR) \citep{NapovPariwise}, the cycle type is the V-cycle, the smoothing operator is PGS-SCM, the degree of freedom of the coarsest space is set to be 10000, and the coarsest space solver is a direct solver. For the restarted GMRES($m$) method, the restarting number $m$ is 30, the maximum number of iterations \emph{MaxIt} is 1000, and the tolerance for relative residual \emph{tol} is $10^{-5}$.

To better evaluate the performance of the proposed methods, we also test the same problems with a commercial simulator (2020 version) for comparison. In the commercial simulator, the default solving method and parameters are used, where the maximum number of iterations \emph{MaxIt} is 1000 and the tolerance for relative residual \emph{tol} is $10^{-5}$. 
Here, we compare the experimental results of the GPU version for commercial and our simulators. The numerical experiments are tested on a machine with Intel Xeon Platinum 8260 CPU (32 cores, 2.40GHz), 128GB DRAM, and NVIDIA Tesla T4 GPU (16GB Memory).

\subsection{The modified SPE1 problem}

The three-phase SPE1 example~\cite{Odeh1981-SPE1} is a benchmark problem for testing ten-year dynamic simulations of immiscible gas flooding (one gas injection well and one oil production well). The initial reservoir state is unsaturated. Given the initial oil-gas and oil-water interface depth, as the reservoir pressure decreases, the gas will gradually dissolve into the oil, and this process will affect the stability of the simulator. The horizontal direction of the oil field is a square with side length of 10000 (ft) and the vertical thickness is 100 (ft). The grid size of the original problem is a $10\times10\times3$ orthogonal grid, and we refine it to get a grid of $80\times80\times40$. Here, we perform numerical tests on the refined grid.

\subsubsection{PGS-MC method}
To investigate the convergence behavior and parallel performance of the proposed PGS-MC method, we compared with the parallel GS and Jacobi methods based on natural ordering, denoted as PGS-NO and PJAC-NO, respectively. The MSP-GMRES method is used as a solver for petroleum reservoir simulation. In the MSP preconditioner, the smoothing operator of the AMG method uses PJAC-NO, PGS-NO, and PGS-MC, denoted as MSP-GMRES-PJAC-NO, MSP-GMRES-PGS-NO, and MSP-GMRES-PGS-MC, respectively. Below, we simulate 10 years (3655.5 days) using the SPE1 example and test the GPU-based parallel performance of the three solvers.

\begin{table}[htpb]
	\centering
	\setlength{\tabcolsep}{4.0pt} 
	\renewcommand\arraystretch{1.5} 
	\caption{Iter and Time (s) of the three solvers for the three-phase SPE1 problem.}
	\label{tab:MC-GS-SPE1}
	\begin{tabular}{c|c|c|c}
		\hline
		Solvers & MSP-GMRES-PJAC-NO  & MSP-GMRES-PGS-NO& MSP-GMRES-PGS-MC \\
		\hline
		Iter &20880 &19339 & 18343 \\
		\hline
		Time &920.83  &888.69  &885.07 \\
		\hline
	\end{tabular}
\end{table}

\reftab{tab:MC-GS-SPE1} presents the total number of linear iterations (Iter) and the total solution time (Time) for the three solvers. We can observe that MSP-GMRES-PJAC-NO has the more number of iterations, our MSP-GMRES-PGS-MC has the least number of iterations (and less time), while MSP-GMRES-PGS-NO's iterations is between MSP-GMRES-PGS-MC and MSP-GMRES-PJAC-NO. This shows that MSP-GMRES-PGS-MC produces the same convergence behavior as the corresponding single-thread algorithm (this conclusion is also confirmed by the OpenMP version~\cite{LiZhao2022}).

\subsubsection{ASMSP-GMRES-GPU method}
\begin{figure}[H]
	\centering	
	\subfigure[Oil production rate]{
		\begin{minipage}[t]{0.45\linewidth}
			\centering
			\includegraphics[width=6cm]{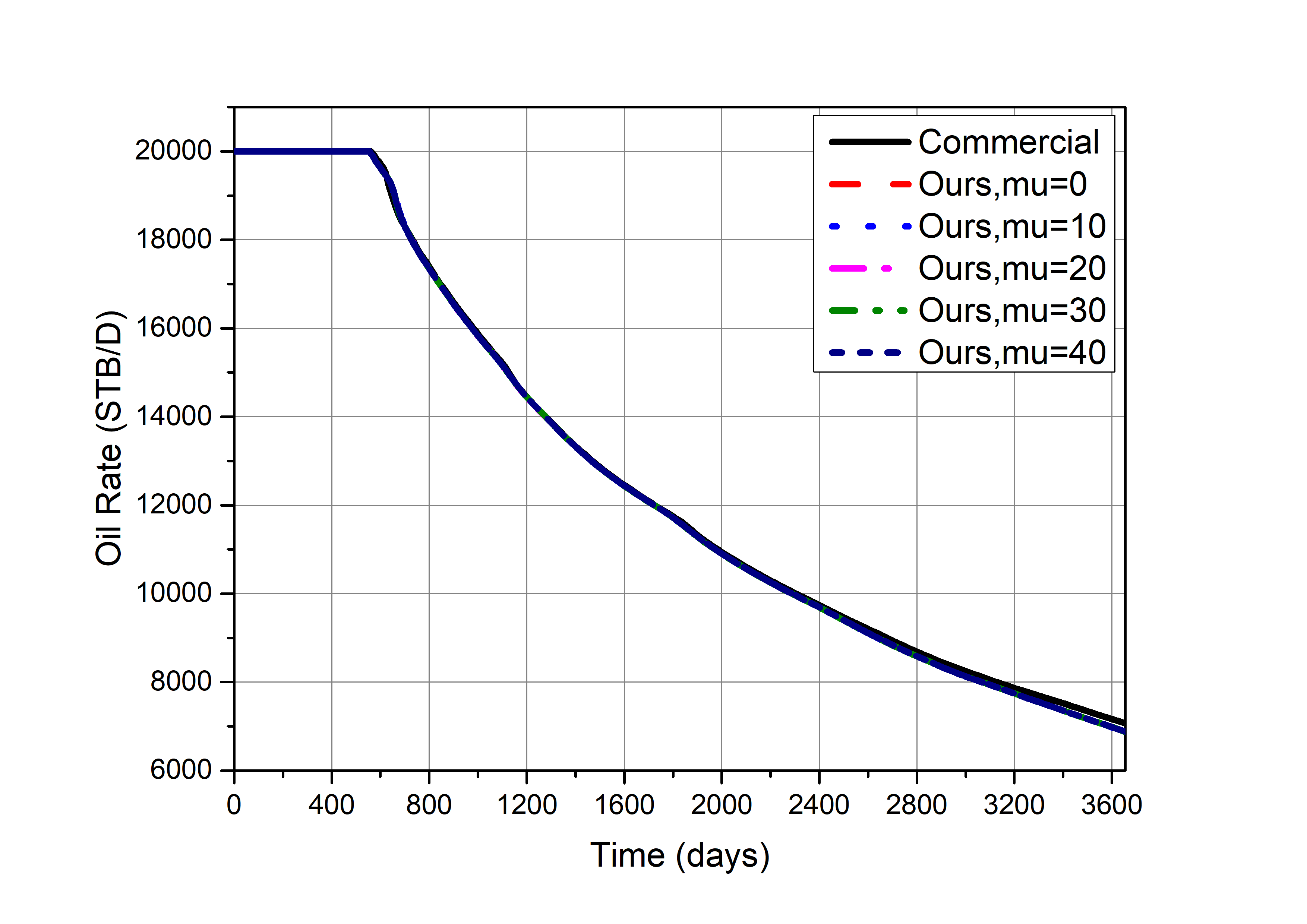}
		\end{minipage}
	}%
	\subfigure[Average pressure]{
		\begin{minipage}[t]{0.45\linewidth}
			\centering
			\includegraphics[width=6cm]{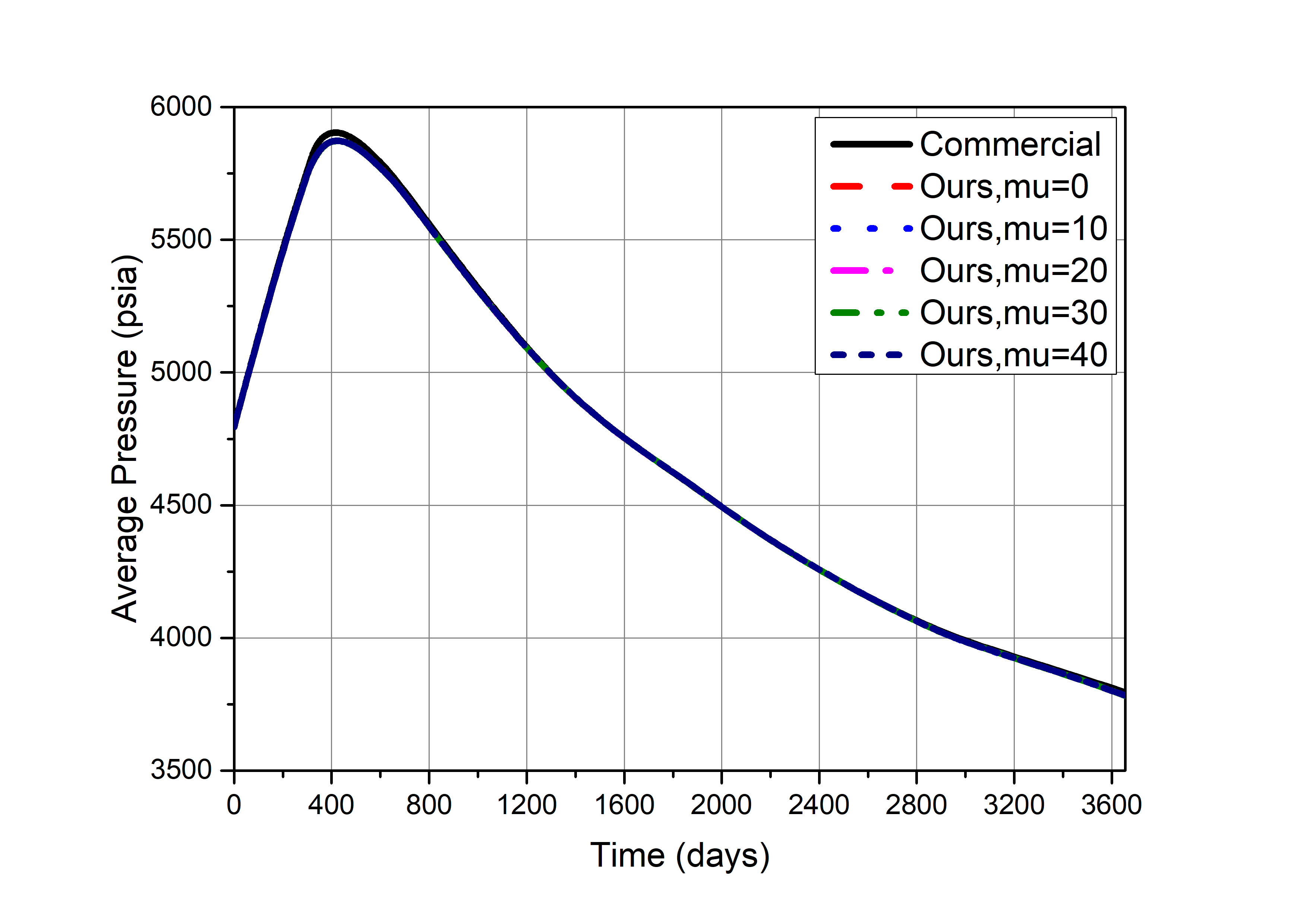}
		\end{minipage}
	}%
	\caption{Field oil production rate and average pressure for the modified SPE1 problem on GPUs.}\label{fig:RES-SPE1-GPU}
\end{figure}

We investigate the effects of $ \mu $ ($ \mu = 0, 10, 20, 30$, and $40$) on the parallel performance of ASMSP-GMRES-GPU. 
We first verify the correctness of ASMSP-GMRES-GPU by comparing our results with commercial simulator. \reffig{fig:RES-SPE1-GPU} shows the field oil production rate and average pressure graphs of five groups $\mu $.

From \reffig{fig:RES-SPE1-GPU}, we can observe that the field oil production rate and average pressure obtained by our and commercial simulators are consistent, indicating that the correctness of our proposed methods is guaranteed.

In addition, to evaluate the parallel performance of the proposed methods for ASMSP-GMRES-GPU, \reftab{tab:GPU-results-SPE1} lists the number of SETUP calls ({SetupCalls}), the ratio of SETUP in the total solution time ({SetupRatio}), the total number of linear iterations ({Iter}), total solution time ({Time}), and parallel speedup ({Speedup}). ASMSP-GMRES-SEQ is the sequential program for reference.
\begin{table}[htpb]
	\centering
	\setlength{\tabcolsep}{4.5pt} 
	\renewcommand\arraystretch{1.5} 
	\caption{SetupCalls, SetupRatio, Iter, Time (s), and Speedup of the different $\mu$ for the SPE1 problem.}
	\label{tab:GPU-results-SPE1}
	\begin{tabular}{c|c|c|c|c|c|c}
	\hline
	Solvers & $\mu$& SetupCalls&	SetupRatio & Iter&	Time&	Speedup \\
	\hline
	ASMSP-GMRES-SEQ &0	  &1178	& 14.48\% &17486	&6094.14 	&---  \\
	\hline
	\multirow{5}{*}{ASMSP-GMRES-CUDA}
	&0		&1177	&61.98\%	&18343	&885.07  	& 6.89 \\
	&10		&1031	&61.26\%	&18529	&867.66 	& 7.02 \\
	&\textbf{20}		&\textbf{236}	&\textbf{50.45\%}	&\textbf{20582}	&\textbf{749.23} 	& \textbf{8.13} \\
	&30		&52		&45.46\%	&22607	&750.45 	& 8.12 \\
	&40		&35		&43.47\%	&24217	&778.65  	& 7.83 \\
	\hline
	\end{tabular}
\end{table}

As can be seen from \reftab{tab:GPU-results-SPE1}, we observe the ASMSP-GMRES-CUDA method. As $\mu$ increases, both the number of SETUP calls and the ratio of SETUP in the total solution time decrease, and the parallel speedup first increases and then decreases (since the number of linear iterations increases gradually). Compared with ASMSP-GMRES-SEQ, when $ \mu=0 $, the speedup of ASMSP-GMRES-CUDA reaches 6.89. In particular, the solution time of ASMSP-GMRES-CUDA with $ \mu=20 $ is reduced from 885.07s to 749.23s compared with $ \mu=0 $. Simultaneously, the speedup of ASMSP-GMRES-CUDA reaches 8.13 compared to ASMSP-GMRES-SEQ. Therefore, the proposed methods can obtain acceleration effects for GPU architecture.

Finally, \reftab{tab:OCP-GPU-SPE1} presents the number of time steps ({NumTSteps}), the number of Newton iterations ({NumNSteps}), the number of linear iterations ({Iter}), the average number of linear iterations per Newton iteration ({AvgIter}), the total simulation time ({Time}), and the parallel speedup ({Speedup}) for the commercial and our simulators, respectively. 
\begin{table}[htpb]
	\centering
	\setlength{\tabcolsep}{4.5pt} 
	\renewcommand\arraystretch{1.5} 
	\caption{NumTSteps, NumNSteps, Iter, AvgIter, Time(s), and Speedup comparisons of the commercial and our simulators for the SPE1 problem.}\label{tab:OCP-GPU-SPE1}	
	\begin{tabular}{c|c|c|c|c|c|c|c}
		\hline
		Simulators &$\mu$ &NumTSteps &NumNSteps &Iter &AvgIter	&Time &Speedup \\
		\hline
		Commercial &---	&410 & 838 &92373 &110.2  &3382.00 &---   \\
		\hline
		\multirow{5}{*}{Ours}
		&0	&267	&1177 &18343 &15.6 &2414.36  &1.40 \\
		&10	&267	&1177 &18529 &15.7 &2344.43 &1.44 \\
		&\textbf{20}	&\textbf{267}	&\textbf{1179}	&\textbf{20582} &\textbf{17.5} &\textbf{2228.08} &\textbf{1.52} \\
		&30	&267  &1178 &22607 &19.2 &2229.42 &1.52 \\
		&40	&267  &1178 &24217 &20.6 &2339.99 &1.45 \\
		\hline
	\end{tabular}
\end{table}

It can be seen from \reftab{tab:OCP-GPU-SPE1} that the commercial simulator requires more  numbers of time steps and linear iterations, as well as more simulation time, compared to our simulator. They yield the average number of linear iterations per Newton iteration of 110.2, about 7 times as ours (when $\mu=0$). 
When $\mu=0$, the speedup of our simulator achieves 1.40 compared to the commercial simulator. When $\mu=20$, the minimum simulation time is 2228.08s, and the speedup is 1.52. This indicates that the proposed methods can improve parallel performance. Finally, it is worth noting that we only parallelize the linear solver in our simulator, while the rest of the simulator is still sequential.

\subsection{The SPE10 problem}
The two-phase SPE10~\citep{2001TenthSPE} benchmark problem with strong heterogeneity is tested to demonstrate the effectiveness of the proposed methods. Its model dimensions are $1200 \times 2200 \times 170$ (ft) and the number of grid cells is $60\times220\times80$ (the total number of grid cells is 1,122,000 and the number of active cells is 1,094,422). In this example, the numerical simulation is carried out for 2000 days. We analyze the effects of $ \mu $ ($ \mu = 0, 10, 20, 30$, and $40$) on the parallel performance of ASMSP-GMRES-GPU. 

\begin{figure}[H]
	\centering	
	\subfigure[Oil production rate]{
		\begin{minipage}[t]{0.45\linewidth}
			\centering
			\includegraphics[width=6cm]{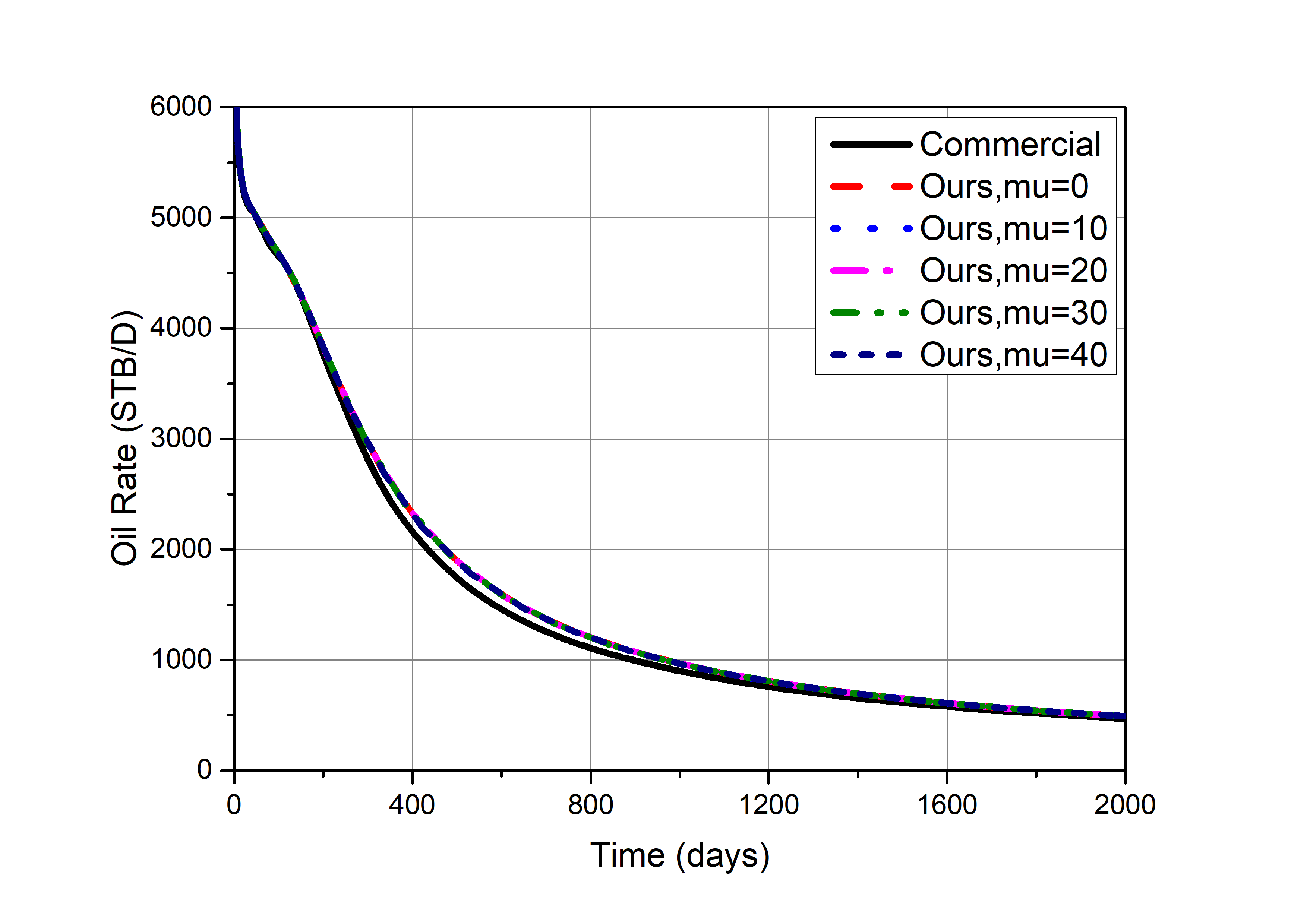}
		\end{minipage}
	}%
	\subfigure[Average pressure]{
		\begin{minipage}[t]{0.45\linewidth}
			\centering
			\includegraphics[width=6cm]{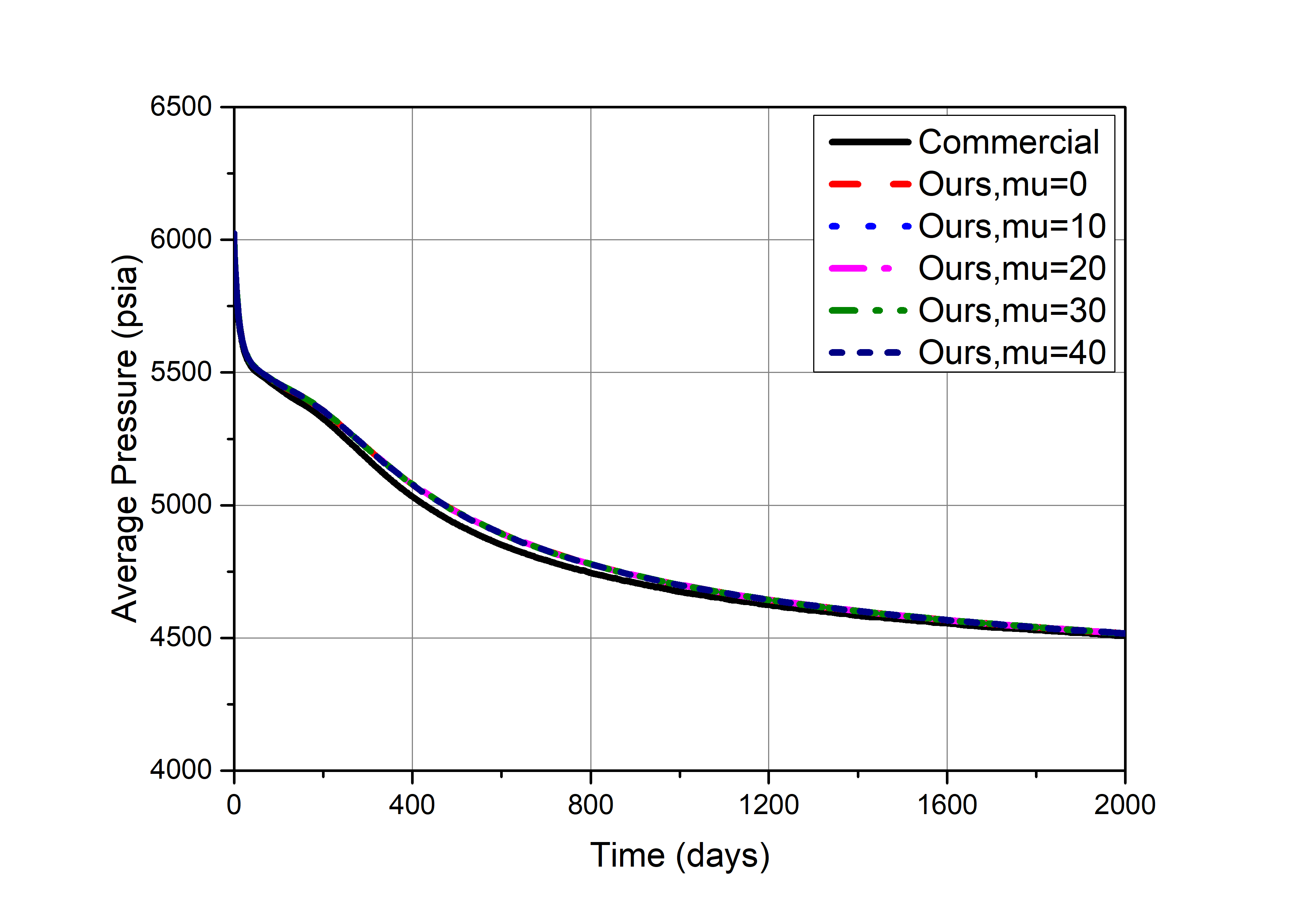}
		\end{minipage}
	}%
	\caption{Field oil production rate and average pressure of the SPE10 problem on GPUs.}\label{fig:FOPR-FPR-SPE10-2P}
\end{figure}

To begin with, we verify the correctness of ASMSP-GMRES-GPU by comparing our results with commercial simulator. The field oil production rate and average pressure graphs of five groups $\mu $ are presented in \reffig{fig:FOPR-FPR-SPE10-2P}. We note that the field oil production rate and average pressure obtained by our and commercial simulators are consistent, indicating that the proposed methods are corrected.

Furthermore, \reftab{tab:GPU-results-SPE10} lists the number of SETUP calls ({SetupCalls}), the ratio of SETUP in the total solution time ({SetupRatio}), the total number of linear iterations ({Iter}), total solution time ({Time}), and parallel speedup ({Speedup}), to assess the parallel performance of the proposed methods. Also, ASMSP-GMRES-SEQ is the sequential program for reference.

\begin{table}[htpb]
	\centering
	\setlength{\tabcolsep}{5.0pt} 
	\renewcommand\arraystretch{1.5}
	\caption{SetupCalls, SetupRatio, Iter, Time (s), and Speedup of the different $\mu$ for the two-phase SPE10 problem.}
	\label{tab:GPU-results-SPE10}
	\begin{tabular}{c|c|c|c|c|c|c}
		\hline
		Solvers & $\mu$& {SetupCalls}&	{SetupRatio} & {Iter}&	{Time}&	{Speedup} \\
		\hline
		ASMSP-GMRES-SEQ &0	&219	&12.40\%	&4252	&4795.04 & ---   \\
		\hline
		\multirow{5}{*}{ASMSP-GMRES-CUDA}
		&0	&219	&55.11\%	&5073	&716.08 	&6.70  \\
		&10	&172	&51.07\%	&5081	&656.37 	&7.31  \\
		&20	&118	&41.36\%	&5863	&630.31 	&7.61  \\
		&\textbf{30}	&\textbf{65}		&\textbf{32.92\%}	&\textbf{6043}	&\textbf{566.68} 	&\textbf{8.46} \\
    &40	&51		&27.87\%	&7055	&610.40 	&7.86 \\
		\hline
	\end{tabular}
\end{table}

From \reftab{tab:GPU-results-SPE10}, we observe the ASMSP-GMRES-CUDA method. As $\mu$ increases, both the number of SETUP calls and the ratio of SETUP in the total solution time decrease, and the parallel speedup first increases and then decreases. Compared with ASMSP-GMRES-SEQ, when $ \mu=0 $, the speedup of ASMSP-GMRES-CUDA reaches 6.70. In particular, the solution time of ASMSP-GMRES-CUDA with $ \mu=30 $ is reduced from 716.08s to 566.68s compared with $ \mu=0 $. Simultaneously, the speedup of ASMSP-GMRES-CUDA reaches 8.46 compared to ASMSP-GMRES-SEQ.

Finally, \reftab{tab:GPU-SPE10-2P} presents the number of time steps ({NumTSteps}), the number of Newton iterations ({NumNSteps}), the number of linear iterations ({Iter}), the average number of linear iterations per Newton iteration ({AvgIter}), the total simulation time ({Time}), and the parallel speedup ({Speedup}) for the commercial and our simulators, respectively. 

\begin{table}[htpb]
	\centering
	\setlength{\tabcolsep}{4.5pt} 
	\renewcommand\arraystretch{1.5} 
	\caption{NumTSteps, NumNSteps, Iter, AvgIter, Time (s), and Speedup comparisons of the commercial and our simulators for the two-phase SPE10 problem.}\label{tab:GPU-SPE10-2P}	
	\begin{tabular}{c|c|c|c|c|c|c|c}
		\hline
		Simulators	&$\mu$ &NumTSteps&NumNSteps&Iter		&AvgIter	&Time &Speedup \\
		\hline
		Commercial	&---	&1004 &1431 &170276	 &119.0 &11034.00 	&---      \\
		\hline
		\multirow{5}{*}{Ours}
		&0	 &53  &219	&5073 &23.2  &1589.62 & 6.94\\
		&10	 &53  &219	&5081 &23.2  &1527.12 & 7.23\\
		&20	 &53  &222	&5863 &26.4  &1534.52 & 7.19\\
		&\textbf{30}	 &\textbf{53}  &\textbf{222}	&\textbf{6043} &\textbf{27.2}  &\textbf{1458.14} &\textbf{7.57} \\
		&40  &54  &221  &7055 &31.9  &1561.90 & 7.06\\
		\hline
	\end{tabular}
\end{table}

According to \reftab{tab:GPU-SPE10-2P}, the commercial simulator requires more numbers of time steps, Newton iterations, and linear iterations compared to our simulator. The average number of linear iterations per Newton iteration is 119.0 (over 5 times as ours when $\mu=0$), and the simulation time is 11034.00s (3.06h). When $\mu=0$, the speedup of our simulator achieves 6.94 compared to the commercial simulator. When $\mu=30$, the minimum simulation time is 1458.14s, and the speedup reaches 7.57. These results show that the parallel performance of the proposed methods outperforms commercial simulators.

\subsection{The modified SPE5 problem}
The SPE5~\citep{Killough1987SPE5} example is a compositional reservoir problem, including six components ($\rm{C}_1$, $\rm{C}_3$, $\rm{C}_6$, $\rm{C}_{10}$, $\rm{C}_{15}$, and $\rm{C}_{20}$), injection well (water alternating gas), and production well. Its reservoir domain is $3500\times3500\times100$ (ft), the original orthogonal grid is $7\times7\times3$, and the simulation period is 20 years. Here, to evaluate the performance of the proposed methods for compositional reservoir, we refine the original grid to obtain a $70\times70\times30$ orthogonal grid, and simulate a period of 2 years. We test the effects of $ \mu $ ($ \mu = 0, 10, 15$, and $20$) on the parallel performance of ASMSP-GMRES-GPU.

\begin{figure}[H]
	\centering	
	\subfigure[Oil production rate]{
		\begin{minipage}[t]{0.45\linewidth}
			\centering
			\includegraphics[width=6cm]{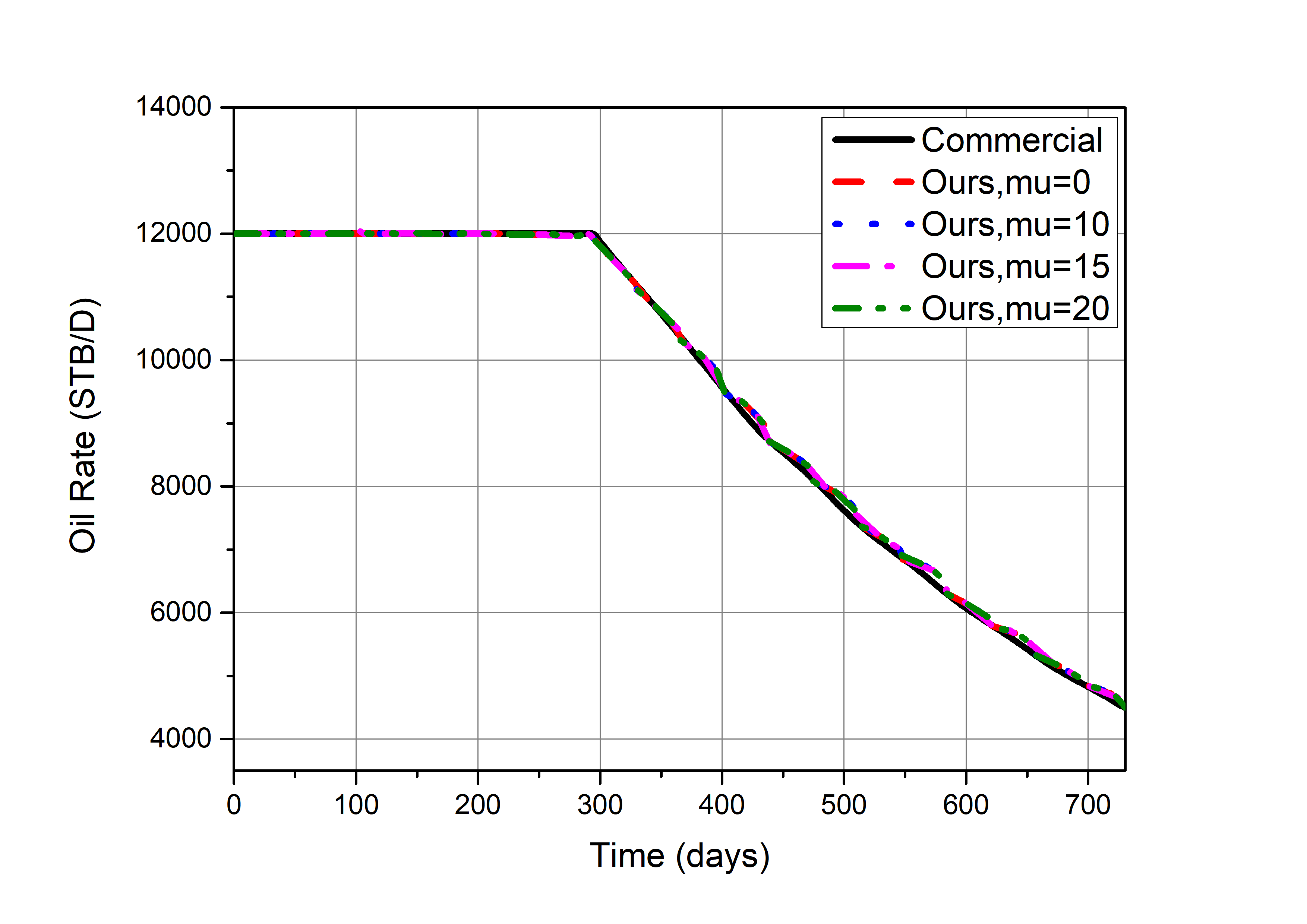}
		\end{minipage}
	}%
	\subfigure[Average pressure]{
		\begin{minipage}[t]{0.45\linewidth}
			\centering
			\includegraphics[width=6cm]{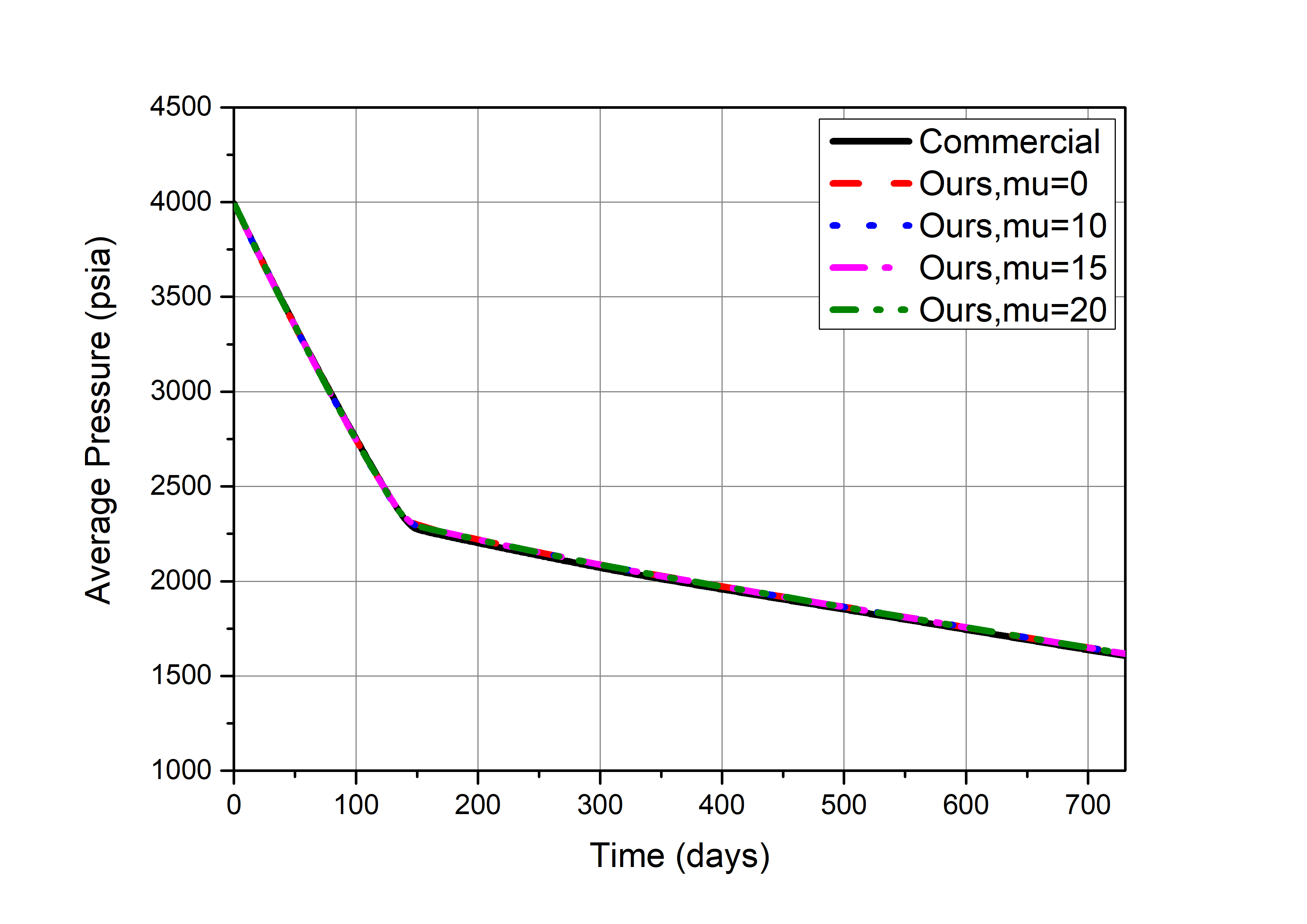}
		\end{minipage}
	}%
	\caption{Field oil production rate and average pressure of the SPE5 problem on GPUs.}\label{fig:FOPR-FPR-SPE5}
\end{figure}

Firstly, we verify the correctness of ASMSP-GMRES-GPU by comparing our results with the commercial simulator. The field oil production rate and average pressure graphs of different $\mu $ are presented in \reffig{fig:FOPR-FPR-SPE5}. We can see the difference in the field oil production rate and average pressure  obtained by our and commercial simulators are consistent, indicating that the proposed methods are corrected.

Furthermore, \reftab{tab:GPU-results-SPE5} lists the number of SETUP calls ({SetupCalls}), the ratio of SETUP in the total solution time ({SetupRatio}), the total number of linear iterations ({Iter}), total solution time ({Time}), and parallel speedup ({Speedup}), to assess the parallel performance of the proposed methods. Also, ASMSP-GMRES-SEQ is the sequential program for reference.

\begin{table}[htpb]
	\centering
	\setlength{\tabcolsep}{5.0pt} 
	\renewcommand\arraystretch{1.5}
	\caption{SetupCalls, SetupRatio, Iter, Time (s), and Speedup of the different $\mu$ for the SPE5 problem.}
	\label{tab:GPU-results-SPE5}
	\begin{tabular}{c|c|c|c|c|c|c}
		\hline
		Solvers & $\mu$& {SetupCalls}&	{SetupRatio} & {Iter}&	{Time}&	{Speedup} \\
		\hline
		ASMSP-GMRES-SEQ &0	&389	&16.30\%	&3747	&	2601.11 & ---   \\
		\hline
		\multirow{4}{*}{ASMSP-GMRES-CUDA}
		&0	&389	&51.43\%	&3969	&341.25 	&7.62  \\
		&10	&186	&34.59\%	&4064	&324.76 	&8.01  \\
		&\textbf{15}	&\textbf{44}		&\textbf{21.72\%}	&\textbf{4508}	&\textbf{313.23} 	&\textbf{8.30}  \\
		&20	&12		&18.01\%	&4747	&314.50 	&8.27 \\
		\hline
	\end{tabular}
\end{table}

From \reftab{tab:GPU-results-SPE5}, we observe the ASMSP-GMRES-CUDA method. As $\mu$ increases, both the number of SETUP calls and the ratio of SETUP in the total solution time decrease, and the parallel speedup first increases and then decreases. Compared with ASMSP-GMRES-SEQ, when $ \mu=0 $, the speedup of ASMSP-GMRES-CUDA reaches 7.62. In particular, the solution time of ASMSP-GMRES-CUDA with $ \mu=15 $ is reduced from 341.25s to 313.23s compared with $ \mu=0 $. Simultaneously, the speedup of ASMSP-GMRES-CUDA reaches 8.30 compared to ASMSP-GMRES-SEQ.

Finally, \reftab{tab:OCP-GPU-SPE5P} presents the number of time steps ({NumTSteps}), the number of Newton iterations ({NumNSteps}), the number of linear iterations ({Iter}), the average number of linear iterations per Newton iteration ({AvgIter}), the total simulation time ({Time}), and the speedup ({Speedup}) for the commercial and our simulators, respectively. 

\begin{table}[htpb]
	\centering
	\setlength{\tabcolsep}{4.5pt} 
	\renewcommand\arraystretch{1.5} 
	\caption{NumTSteps, NumNSteps, Iter, AvgIter, Time (s), and Speedup comparisons of the commercial and our simulators for the SPE5 problem.}\label{tab:OCP-GPU-SPE5P}	
	\begin{tabular}{c|c|c|c|c|c|c|c}
		\hline
		Simulators	&$\mu$ &NumTSteps&NumNSteps&Iter		&AvgIter	&Time &Speedup \\
		\hline
		Commercial	&---	&382 &748 &47027	 &62.9  &2339.00 	&---      \\
		\hline
		\multirow{4}{*}{Ours}
		&0	 &147  &389	&3969 &10.2  &2178.78 & 1.07 \\
		&10	 &147  &389	&4064 &10.4  &2159.22 & 1.08 \\
		&\textbf{15}	 &\textbf{147}  &\textbf{389}	&\textbf{4508} &\textbf{11.6}  &\textbf{2142.22} & \textbf{1.09} \\
		&20  &147  &389 &4747 &12.2  &2143.47 & 1.09 \\
		\hline
	\end{tabular}
\end{table}

According to \reftab{tab:OCP-GPU-SPE5P}, the commercial simulator requires more numbers of time steps, Newton iterations, and linear iterations compared to our simulator. The average number of linear iterations per Newton iteration is 62.9 (over 6 times as ours when $\mu=0$), and the simulation time is 2339.00s. When $\mu=0$, the speedup of our simulator achieves 1.07 compared to the commercial simulator. When $\mu=15$, the minimum simulation time is 2142.22s, and the speedup reaches 1.09. 
Finally, we discuss the reasons for the low speedup in this example. As the number of components increases in the compositional model, the complexity of the kernel algorithms increases. Moreover, we only parallelize the linear solver in our simulator, while the rest of the simulator is still sequential. As a result, the parallel solver time is only about 15\% of the total simulation time in this example. Hence the linear solver part is not the main computational bottleneck.

\section{Conclusions}\label{sec:6}
In this work, we developed a parallel multi-stage preconditioner for the system of linear algebraic equations arising from the fully implicit approach for the compositional model. We proposed an efficient multi-stage preconditioner with an adaptive SETUP to improve the parallel performance of the preconditioner. Furthermore, we developed an improved parallel GS algorithm based on the adjacency matrix. This algorithm can be applied to the smoothing operator of the AMG methods and yielded the same convergence behavior as the corresponding single-threaded algorithm. We believe the proposed parallel solver framework will provide a feasible approach to the efficient numerical solution of various application problems. In the future, we will further improve the proposed solver and parallelize our compositional simulator OpenCAEPoro.

\backmatter

%\bmhead{Supplementary information}
%
%If your article has accompanying supplementary file/s please state so here. 
%
%Authors reporting data from electrophoretic gels and blots should supply the full unprocessed scans for key as part of their Supplementary information. This may be requested by the editorial team/s if it is missing.
%
%Please refer to Journal-level guidance for any specific requirements.

\bmhead{Acknowledgments}

This work was supported by the Postgraduate Scientific Research Innovation Project of Hunan Province (No. CX20210607) and Postgraduate Scientific Research Innovation Project of Xiangtan University (No. XDCX2021B110). Zhang was partially supported by the National Science Foundation of China (No. 11971472). Feng was partially supported by the Excellent Youth Foundation of SINOPEC (No. P20009). Shu was partially supported by the National Science Foundation of China (No. 11971414).

\bmhead{Statements and Declarations}
On behalf of all authors, the corresponding author states that there is no conflict of interest.

\bibliography{references}


\end{document}